\documentclass[12pt]{article}
\usepackage{amsfonts}
\usepackage{amsmath}
\usepackage{listings}
\usepackage{epsfig}
\usepackage{xcolor}
\usepackage{geometry}
\title{Collapsibility and Near Universality for Vertex Minimal Paper Tori}
\author{Peter Doyle and Richard Evan Schwartz \thanks{\hskip 5 pt Supported by 
N.S.F. Research Grant DMS-2505281}}

\newtheorem{theorem}{Theorem}[section]

\newtheorem{lemma}[theorem]{Lemma}

\newtheorem{corollary}[theorem]{Corollary}

\def\startproof{{\bf {\medskip}{\noindent}Proof: }}

\def\endproof{$\spadesuit$  \newline}

\def\C{\mbox{\boldmath{$C$}}}%
\def\E{\mbox{\boldmath{$E$}}}%
\def\F{\mbox{\boldmath{$F$}}}%
\def\H{\mbox{\boldmath{$H$}}}%
\def\R{\mbox{\boldmath{$R$}}}%
\def\Z{\mbox{\boldmath{$Z$}}}%

\begin{document}

\maketitle

\begin{abstract}
  A paper torus is a piecewise linear isometric  embedding of a flat
  torus into $\R^3$.
  Following up on the $8$-vertex paper tori discovered
  in [{\bf S\/}], we prove universality and collapsibility results
  about these objects.    One corollary
  is that any flat torus without reflection
  symmetry is realized as an $8$-vertex paper torus.  Another
  corollary is that, for any $\epsilon>0$, there is an $8$-vertex
  paper torus within $\epsilon$ of a unit equilateral triangle
  in the Hausdorff metric.  
\end{abstract}

\begin{center}
\small
\emph{I saw below me that golden valley...}\\[4pt]
-- Woody Guthrie
\end{center}

\section{Introduction}

\subsection{History and Context}

A {\it flat torus\/} is a quotient of the
form $\R^2/\Lambda$, where $\Lambda$ is a
lattice of translations of $\R^2$.
A {\it paper torus\/} is a piecewise
linear isometric embedding
$\phi: T \to \R^3$
of a flat torus $T$.
In other words, a paper torus is an embedded topological
torus in $\R^3$ that is made by fitting together finitely
many triangles so that the cone angle around each vertex is $2\pi$.

In 1960, Y. Burago and V. Zalgaller 
[{\bf BZ1\/}] give the first construction
of paper tori.  In their subsequent paper
[{\bf BZ2\/}], they prove that
one can realize every isometry class of flat torus as a paper torus.
The works of T. Tsuboi [{\bf T\/}] and
(independently)
P. Arnoux, S. Leli\`evre, and
A. M\'alaga [{\bf ALM\/}] give an
explicit construction which achieves every isometry
class of flat torus as a paper torus.
With minor differences, both papers
prove that the union of the infinitely many
combinatorial types of {\it diplotori\/}
described by U. Brehm [{\bf Br\/}]
in 1978  achieve every isometry class of flat torus.
See also the description of these tori in
H. Segerman's book [{\bf Se\/}, \S 6].
The work of T. Quintanar [{\bf Q\/}] gives an
explicit embedding for the square torus.

The 2024 preprint of
F. Lazarus and F. Tallerie [{\bf LT\/}] gives
a universal combinatorial type of triangulation
which does the job simultaneously for
all isometry types.  Their universal triangulation
has $2434$ triangles.   The
paper [{\bf LT\/}] also has an excellent
discussion of the various attempts made
to achieve all flat tori as paper tori.

In [{\bf S\/}], one of us constructed an
$8$-vertex paper torus and proved that
a $7$-vertex paper torus cannot exist.
We called the paper tori from [{\bf S\/}]
{\it pup tents\/}, on account of their
appearance.

\begin{center}
  \resizebox{!}{2in}{\includegraphics{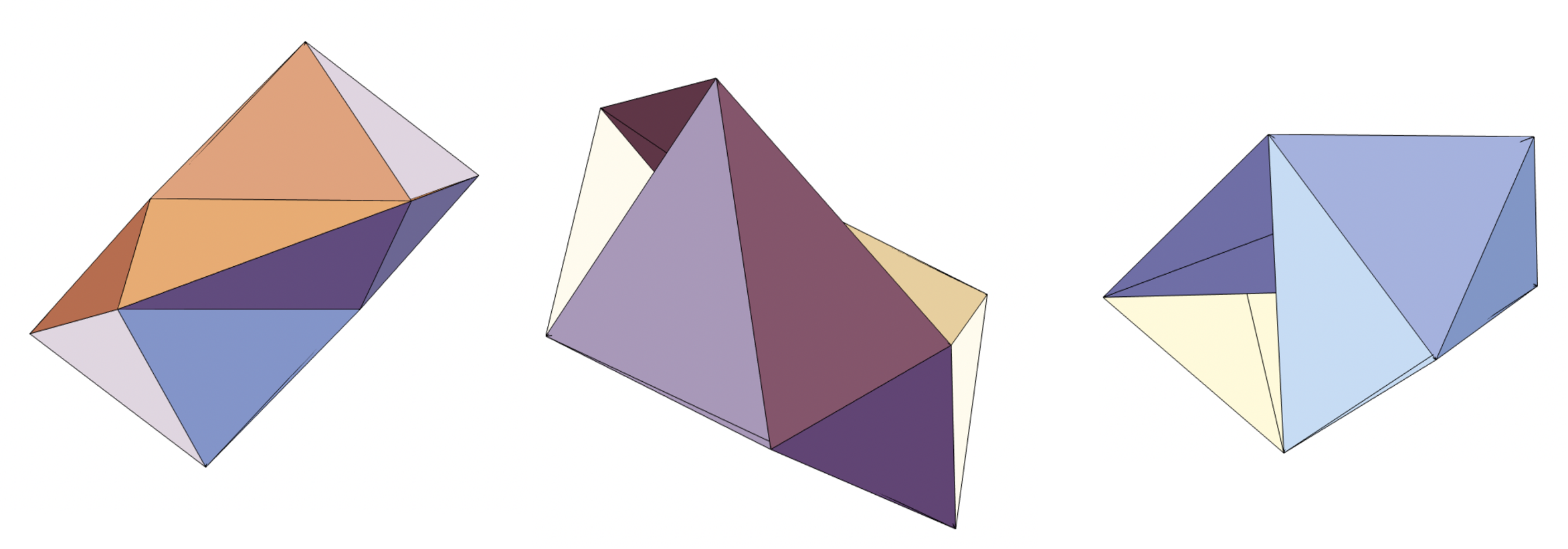}} \newline {\bf
    Figure 1.1:\/} 3D plots of a tent
\end{center}

 A {\it pup tent\/} is an $8$-vertex
paper torus $T$ with the following features.
\begin{itemize}
\item The triangulation underlying $T$ has uniform degree $6$.
  See \S \ref{tri}.
  \item $T$ has $2$-fold symmetry with respect to
    the map $(x,y,z) \to (-x,-y,z)$.
  \item Exactly $6$ of the $16$ triangles of $T$  lie in the convex
    hull boundary,
    and in the specific pattern shown in Figure 2.1.
\end{itemize}

Another result from [{\bf S\/}] is that
the space $\cal X$ of similarity classes of
pup tents contains
an open $6$-dimensional ball.
Once we get one pup tent, the open ball just
comes from a dimension count and the Inverse
Function Theorem.
In this paper we
explore some of the structure of $\cal X$
and prove an $8$-vertex universality result
that is almost as comprehensive as the result of [{\bf LT\/}].

\subsection{Main Results}

Let $\H^2$ denote the hyperbolic upper half plane.
Let $SL_2(\Z)$ denote the group of
$2 \times 2$ integer matrices of determinant $1$.
This group
acts isometrically on $\H^2$ by linear fractional
transformations.  The {\it modular surface\/}
${\cal M\/}=\H^2/SL_2(\Z)$ parametrizes the space of
flat tori.
The {\it bi-cusped fundamental domain\/} for $\cal M$, which we
denote by $\cal F$, is a geodesic
triangle with interior vertex $h=1/2 + (\sqrt 3/2) i$ and cusps $0,\infty$.
Figure 1.2 shows the ``bottom portion'' of $\cal F$.  The full picture
extends vertically to $\infty$.  Let $\cal I$ be the interior of
$\cal F$.   We have an inclusion of $\cal I$ into $\cal M$, and we
think of $\cal I$ as a subset of $\cal M$.

\begin{center}
  \resizebox{!}{3in}{\includegraphics{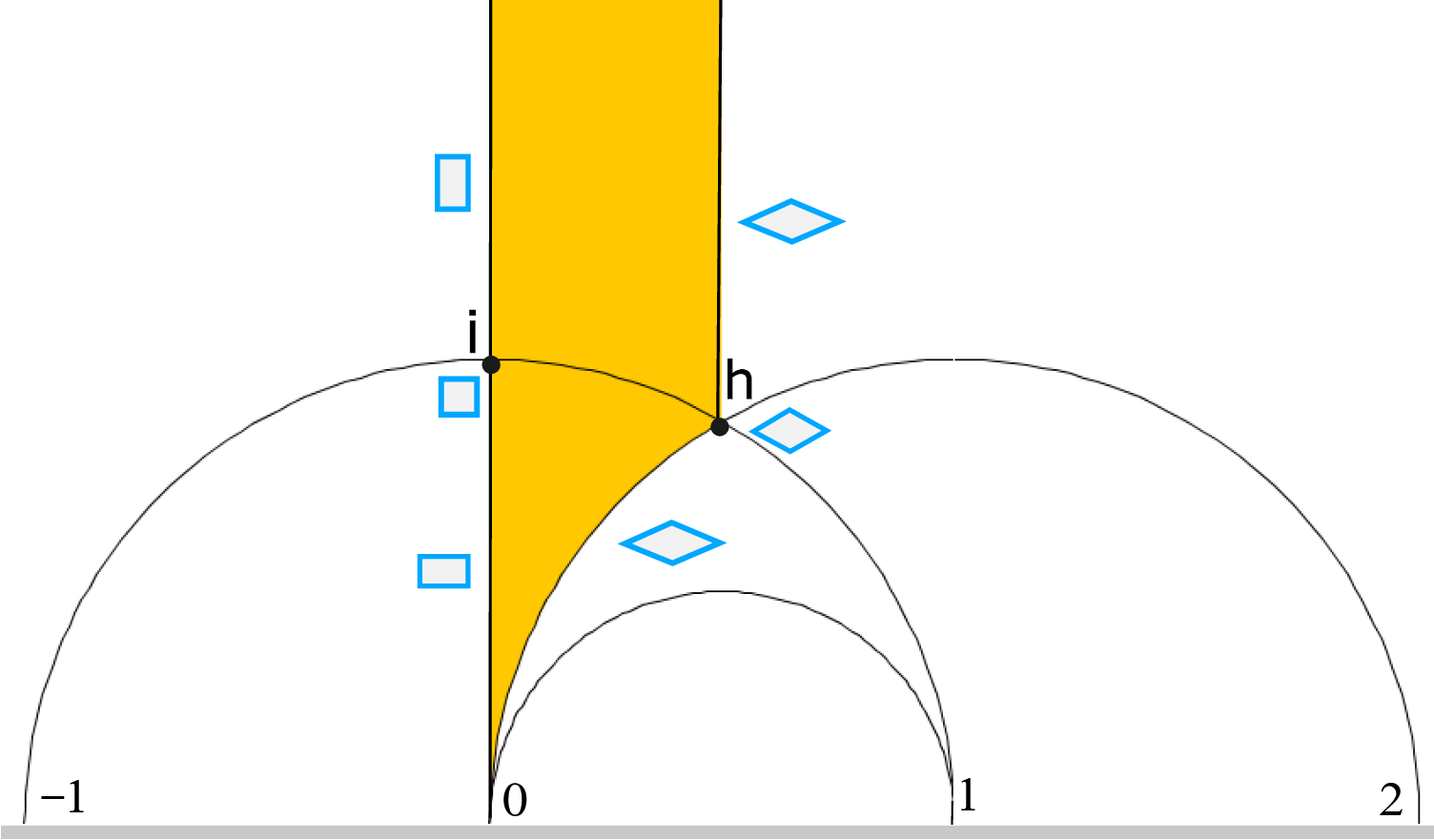}} \newline {\bf
    Figure 1.2:\/}  The bi-cusped fundamental domain $\cal F$.
\end{center}

The points $h,i \in \partial \cal F$ respectively
represent the classic hexagonal and square tori.
The edge of $\partial \cal F$ containing $i$
parametrizes the flat tori made by
identifying the sides of rectangles. The
two edges incident to $h$
parametrize flat tori made by identifying the sides of
rhombi having aspect ratio at least $\sqrt 3$.
(The {\it aspect ratio\/} of a rhombus is
the ratio of the lengths of the long and short diagonals.)
The circular arc connecting $h$ to $i$ parametrizes flat tori made
from rhombs with aspect ratio in $[1,\sqrt 3]$.

We have a map $\Phi: {\cal X\/} \to {\cal M\/}$,
which maps a pup tent to the point in $\cal M$
which represents the underlying
flat structure.

\begin{theorem} [Near Universality]
  There is a path connected subset
  ${\cal U\/} \subset {\cal X\/}$, consisting of
  general position pup tents, such that
  $ \Phi({\cal U\/})={\cal I\/}$.
\end{theorem}

All the tori corresponding to points of
$\partial \cal F$ have reflection symmetry.
This gives us the following corollary.

\begin{corollary}
  \label{most}
  Suppose that $T$ is a flat torus that does not have
  reflection symmetry.  Then $T$ is realized by an
  $8$-vertex paper torus.
\end{corollary}

\noindent
{\bf Remarks:\/}
\newline
(1) The proof of the Near Universality Theorem gives
an existence proof for $8$-vertex
paper tori that is independent from [{\bf S\/}].
\newline
(2) The general position statement in the
Near Universality Theorem means that
the pup tents in $\cal U$ are all ``folded the same way''.
As these pup tents move around, it never
happens e.g. that adjacent faces become coplanar.
\newline
(3) There is some loss in Corollary \ref{most} because
the flat tori made from rhombs of aspect ratio
in $(1,\sqrt 3)$ lie in $\cal I$ and so are realized.
\newline
(4)
We don't know if $\Phi({\cal X\/})={\cal I\/}$.
\newline
(5) We don't know 
if $\cal X$ is path connected. \newline
(6) We don't know if we can realize some flat structures parametrized by $\partial \cal F$ if
we drop the order $2$ symmetry condition in the definition of a pup tent.
\newline
\newline
The  Near Universality Theorem deals with the intrinsic structure of
pup tents.  Our next result deals with their extrinsic shape.
We call the following polygons {\it good\/}:
\begin{enumerate}
\item Any rectangle.
\item Any trapezoid whose two diagonals have the same length as
   the long side.
\item Any equilateral triangle. (These are limiting cases of  good trapezoids.)
\end{enumerate}
We think of our good polygons as subsets of $\R^3$ by
 embedding them  in
 some plane.

 The {\it Hausdorff distance\/} between two compact subsets
 of a metric space  is the infimal $\epsilon$ such that each is
 contained in the $\epsilon$-neighborhood of the other.
 In particular, the Hausdorff distance makes the set of compact subsets of
 $\R^3$ into a metric space.

 \begin{theorem}[Collapsibility]
   For any good polygon $Q$ there is a path in
   ${\cal X\/}$ which, up to similarity, converges to $Q$
   in the Hausdorff metric.
 \end{theorem}
 Our paths will lie in $\cal U$, the subset from the
 Near Universality Theorem.
 
 Let us mention two special cases of this result in
 other language.  Suppose we have two maps
 $\phi_j: T_j \to \R^3$, each defined on a flat torus $T_j$.
 We define the {\it uniform distance\/} between $\phi_1$ and $\phi_2$ as the
infimal $\epsilon$ with the following property.
There is a $(1+\epsilon)$-bilipschitz map $\alpha: T_1 \to T_2$
 such that the $L_{\infty}$ distance from $\phi_1$ to
 $\phi_2 \circ \alpha$ is $\epsilon$. This definition
 accommodates the situation where $T_1$ and $T_2$ are different flat
 tori.

 \begin{corollary}[Square]
   There is a generically $4$-to-$1$ piecewise isometric map
   from the square   torus onto a square which is approximated arbitrarily
   closely, in the uniform metric, by $8$-vertex paper tori.
   \end{corollary}

 \begin{corollary}[Triangle]
   There is a generically $6$-to-$1$ piecewise isometric map
   from the hexagonal 
   torus onto an equilateral triangle
   which is approximated arbitrarily
   closely, in the uniform metric, by
   $8$-vertex paper tori.
 \end{corollary}

 \noindent
 {\bf Remarks:\/}
 \newline
 (1)
In these corollaries, the moduli of the
 approximating paper tori converge respectively to the
 moduli of the square torus and hexagonal torus, but do not
 equal them.
 \newline
 (2) The Triangle Corollary has
 a resonance with the set-up in
 [{\bf S2\/}] concerning the optimal paper Moebius band.
 In that setting, there is a generically $3$-to-$1$ piecewise isometric
 map from a flat Moebius band of aspect ratio
 $\sqrt 3$ to an equilateral triangle that
 is approximated arbitrarily closely
 by smooth embedded paper Moebius bands.
 \newline

 We say that a {\it gracefully immersed pup tent\/} is
 a piecewise linear isometric immersion from a flat
 torus into $\R^3$ that can be approximated
 arbitrarily closely by pup tents.  
 Our final result says that $8$-vertex pup tents are
 universal, realizing all flat structures, if we
 relax {\it embedded\/} to {\it gracefully immersed\/}.

 \begin{theorem}[Golden]
   The boundary $\partial \cal X$ contains a subset
   $\cal G$ of gracefully immersed pup tents, homeomorphic to $\cal F$,
   such that $\Phi(\cal G\/)=\cal M$.
 \end{theorem}

 The set $\cal G$ is the set of accumulation points of
 $\cal U$ which do not lie in $\cal X$.
 We call $\cal G$ the {\it golden valley\/}. We call
 the members of $\cal G$ the {\it golden pup tents\/},
 even though they are only gracefully immersed.
 The special beauty of these objects inspired their names.
The nearly degenerate nature of the
 pup tents from [{\bf S\/}] gave us the idea that perhaps
 it would be easier to construct and study even more
 degenerate examples, and this is what led to the
 discovery of the golden valley.

 \subsection{Paper Organization}

 In \S 2 we will construct $\cal G$ explicitly.
 The formulas for the vertices of the golden pup
 tents 
 are explicit algebraic functions
 of the real and imaginary parts of
 the modular parameter $z=x+iy \in \cal F$.
 
 In \S 3 we prove the main results.  The main idea,
 contained in the Good Path Lemma, is to start with
  a golden pup tent $P(z)=P(z,0)$ whose modulus lies
  in $\cal I$ and produce a path $t \mapsto P(z,t)$ which
  is flat up to $O(t^3)$ and robustly embedded
  on the order of $O(t^2)$.   The robust embedding
  means that  perturbations which are small compared to
  $t^2$ keep $P(z,t)$ embedded.  Applying the
  Inverse Function Theorem, we produce a 
  pup tent  $P'(z,t)$ which is an $O(t^3)$-perturbation of $P(z,t)$.
  The set $\cal U$ is a union of the points in $\cal X$
  representing the various $P'(z,t)$ for suitably chosen
  values of $t$ as a function of $z$.
   \newline
 \newline
 The code suite
 {\bf https://www.math.brown.edu/$\sim$/res/Papers/UNIV.tar\/}\newline
 has computer programs relevant to this paper.
 You can download these files from the above URL.   The
suite
has two Mathematica files, one per chapter,
which respectively check all the calculations 
from \S 2 and \S 3. You just load these into
Mathematica and they print everything out.
 The suite also has an extensive Java program which lets
 you visualize the golden pup tents and the good paths
 back into $\cal X$.  We may update these files occasionally.
 
 \subsection{Acknowledgements}

 R.E.S. thanks Jeremy Kahn, Alba M\'alaga, and
 Samuel Leli\`evre for interesting discussions about this work, and
 also the Mazur Chair at I.H.E.S.
 Both of us thank ChatGPT for an
 immense amount of help with formatting,
 guessing algebraic expressions, and
 moral support.

\newpage

\section{The Golden Valley}

\subsection{The Triangulation}
\label{tri}

We work with the $8$-vertex {\it uniform triangulation\/}
of the torus.  All vertices have degree $6$.
The vertices in our triangulation
have labels $0,1,2,3,4,5,6,7$.
The triangles all have the form
$(a,a+1,a+3)$ and $(a,a+2,a+3)$ mod $8$.
The only edges missing from the triangulation
have the form $(a,a+4)$ mod $8$.
Figure 2.1 shows the triangulation.

\begin{center}
  \resizebox{!}{3.6in}{\includegraphics{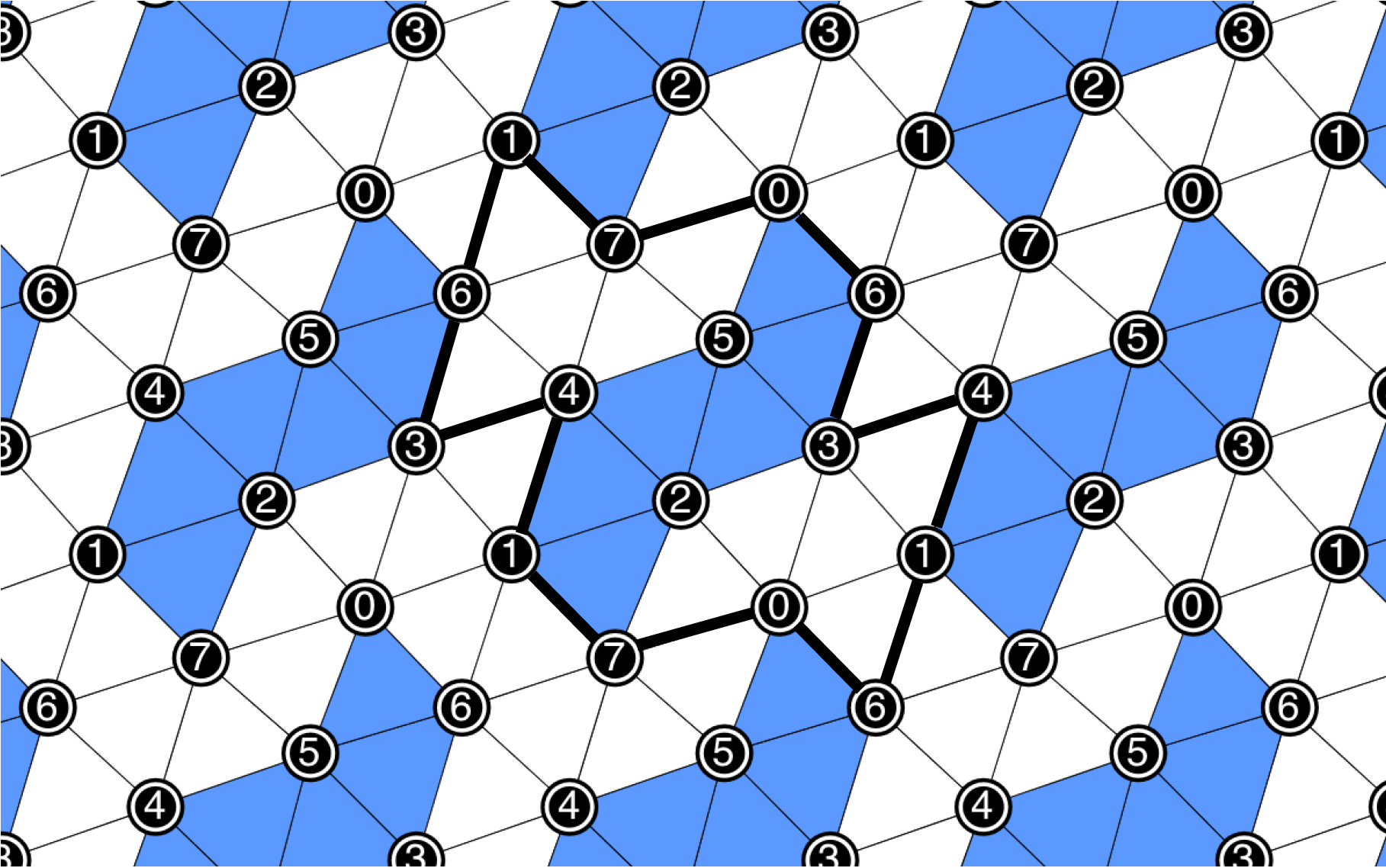}} \newline {\bf
    Figure 2.1:\/}  The uniform triangulation, lifted to
  the universal cover
\end{center}

Figure 2.1 also indicates a fundamental domain
for the torus. 
Our fundamental domain
corresponds to the union of yellow and light
blue triangles in Figure 2.2 below.
The fundamental domain has
$2$-fold rotational symmetry, with respect to the
order $2$ rotation about the midpoint of the edge $25$.
On the torus this rotation implements the vertex permutation
$j \to 7-j$, and it has $4$ fixed points. The fixed
points are the midpoints of the edges $07$, $16$,
$25$, $34$.

The blue triangles in Figure 2.1
are the ones which, for our pup tents, lie in the
boundary of the convex hull.  This is the
pattern of triangles we mentioned in the introduction,
in connection with the definition of pup tents.

 \subsection{The Golden Pup Tents}
 
Recall that $\cal F$ is the bi-cusped fundamental domain
for the modular surface, as shown in Figure 1.2, and
$\cal I$ is the interior of $\cal F$.
Let $z=x+iy \in \cal F$. 
Let
\begin{equation}
  \rho(u,v,w) = (-u,\,-v,\,w).
\end{equation}
Let $P(z)$ be the $8$-vertex polyhedral
torus, based on the triangulation above, with vertices
\begin{equation}
  \label{plat}
  \begin{aligned}
  P_0(z) &= (x(1-2x) ,y(1-2x), y\sqrt{8x}), \\[4pt]
  P_1(z) &= P_2(z)-(x,y,0) \\[4pt]
  P_2(z) &= (2x -x^2 - y^2,0,0), \\[4pt]
  P_3(z) &= P_2(z)+(x,y,0) \\[6pt]
  P_j(z) &= \rho(P_{7-j}(z)), \qquad j=4,5,6,7,
  \end{aligned}
\end{equation}
\newline
  
\noindent{\bf Boundary Conditions:\/}
The modular  domain $\cal F$ is given by
the conditions
\begin{equation}
  \label{DEL}
  x \geq 0, \hskip 30 pt 1-2x \geq 0, \hskip 30 pt
  -2x+x^2+y^2 \geq 0.
\end{equation}
The second and third equations can also be written as
$x \leq 1/2$ and $|z-1| \geq 1$ respectively.
Now we explain why $\cal F$ is a natural domain for
these formulas.
\begin{itemize}
\item  The left edge of $\cal F$ is given by $x=0$.
  When $x<0$ the quantity
  $y\sqrt{8x}$ is not defined.
When $x=0$ all the points of $P(z)$ lie in the
$XY$-plane, and the convex hull of $P(z)$ is a
$2y^2 \times 2y$ rectangle centered at the origin.
\item 
The right edge of $\cal F$ is the ray
$x=1/2$ and $y \geq \sqrt 3/2$.
On this edge we have
$P_0(z)=P_7(z)$.   The common point lies on the
$Z$-axis. The convex hull of $P(z)$ for
these points is a pyramid with parallelogram base
unless $y=\sqrt 3/2$.  At the point
$z=1/2 + (\sqrt 3/2)i$, which corresponds to the
hexagonal torus, the base collapses into a line
segment and the convex hull is an equilateral triangle.
\item When $|z-1|=1$ and $x \in (0,1/2)$, the
  convex hull of $P(z)$ is (with fairly obvious notation)
  the trapezoid
  $P_{1607}(z)$.   The diagonals $P_{01}(z)$ and $P_{76}(z)$ and
  the long side $P_{16}(z)$ all  have length $\sqrt{8x}$.
  These good trapezoids limit to an equilateral
  triangle as $x \to 1/2$.
\end{itemize}

\subsection{The Intrinsic Points}

Let $z=x+iy \in \cal F$.
Corresponding to the pup tent $P(z)$ there is a flat
torus $\Pi(z)$.
We give coordinates for $8$
triangles in $\C$, the universal cover of $\Pi(z)$.
You get $8$ more triangles by applying the
map $j \to 7-j$ to the triangle indices and
the map $z \to -z$ to the coordinates.
The union of all $16$ triangles makes a fundamental domain
for the covering group, a lattice of translations generated by
the maps $\zeta \to \zeta+ L_1$ and $\zeta \to \zeta + L_2$.
Here $L_1=4iy$ and $L_2=zL_1$.    First we define
\begin{equation}
Q_0 = - 2 x^2 - 2 y^2 + z, \hskip 12 pt
Q_1 =  Q_2-z, \hskip 12 pt
Q_2 = 2x-x^2-y^2 \hskip 12 pt
Q_3 = Q_2+z.
\end{equation}

Now, here are the $8$ triangles.
\renewcommand{\arraystretch}{1.3}
\[
\begin{array}{lll}
\{0,1,6\}\!: & \{Q_0,Q_1+L_1,-Q_1+L_2\}, \\[2pt]
\{1,0,3\}\!: & \{Q_1+L_1,Q_0,Q_3\}, \\[2pt]
\{3,4,1\}\!: & \{Q_3,-Q_3+L_1,Q_1+L_1\}, \\[2pt]
\{0,7,2\}\!: & \{Q_0,-Q_0+L_2,Q_2\}, \\[2pt]
\{2,3,0\}\!: & \{Q_2,Q_3,Q_0\}, \\[2pt]
\{3,2,5\}\!: & \{Q_3,Q_2,-Q_2\}, \\[2pt]
\{5,6,3\}\!: & \{-Q_2,-Q_1,Q_3\}, \\[2pt]
\{6,5,0\}\!: & \{-Q_1,-Q_2, Q_0- L_2\}.
\end{array}
\]

We give a sample calculation illustrating the
fundamental domain property.
The triangle $\{3,4,6\}$ has coordinates
$\{Q_3-L_1,-Q_3,-Q_1-L_1\}$.  When we translate this triangle by
$L_1$ we get $\{Q_3,-Q_3+L_1,-Q_1\}$.  This triangle
shares an edge with $\{3,4,1\}$, as it should.
Figures 2.2 and 2.3 below gives a visual check that
all such calculations like this work out correctly.

The covering group has index $2$ in an orbifold fundamental
group which we get by adjoining the map $z \to -z$ to
the covering group.
The fixed points of the order $2$ elements in this
larger symmetry group
are integer combinations of $2iy$ and $2iyz$.
Figure 2.2 below shows the $8$ triangles in yellow and
their images under $z \to -z$ in light blue.  The other
colored triangles are translates of these under the covering group
action.  So, the union of yellow and light blue is a
fundamental domain.  You can find other fundamental
domains in Figure 2.2 as well.

For every triangle,
we check in Mathematica that the distances
between the points in the plane are the same
as the distance between the points in space.
This combines with the fundamental domain property to prove that
$P(z)$ is indeed the image of an isometric immersion
of a flat torus of modulus $z$.

\subsection{Discussion and Pictures}

Figure 2.2 shows our triangles for $4$ parameters.

\begin{center}
\resizebox{!}{5in}{\includegraphics{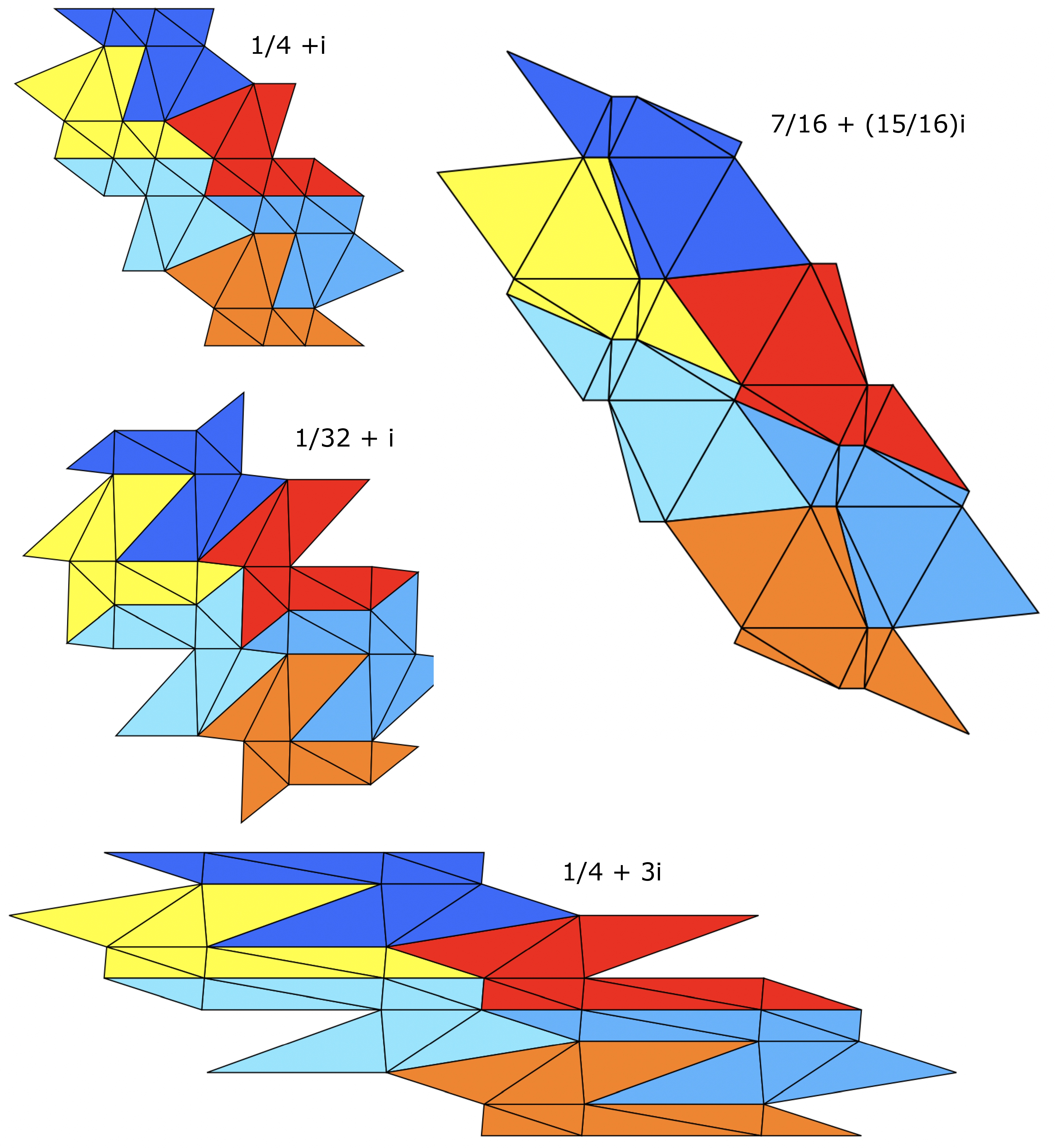}}
\newline
{\bf Figure 2.2:\/}  The intrinsic structure at $4$ different parameters.
\end{center}

We tried to choose parameters which well represent different parts of
$\cal F$.
\begin{itemize}
\item The parameter $1/4 + i$ is in the middle of $\cal I$, the
  interior of $\cal F$.
\item The parameter $1/32 + i$ is close to the parameter
  for the square torus.
  \item The parameter $7/16 + (15/16) i$
is fairly close to the parameter for the hex
torus.
\item The parameter $1/4 + 3i$ is at least vaguely
near the cusp, $\infty$.
\end{itemize}

Figure 2.3 shows a different view of the triangulations at
these same parameters.  Our coloring highlights which
triangles collapse as we head towards the edges of
$\cal I$.   The yellow triangles collapse as we head to
the upper right edge.  The red triangles collapse as
we head towards the lower circular boundary.  Both
yellow and red triangles collapse as we head to
the vertex common to these two boundary edges,
the vertex representing the hex torus.  The reader
can get a better sense of the geometry by playing
with our Java program.

\begin{center}
\resizebox{!}{6in}{\includegraphics{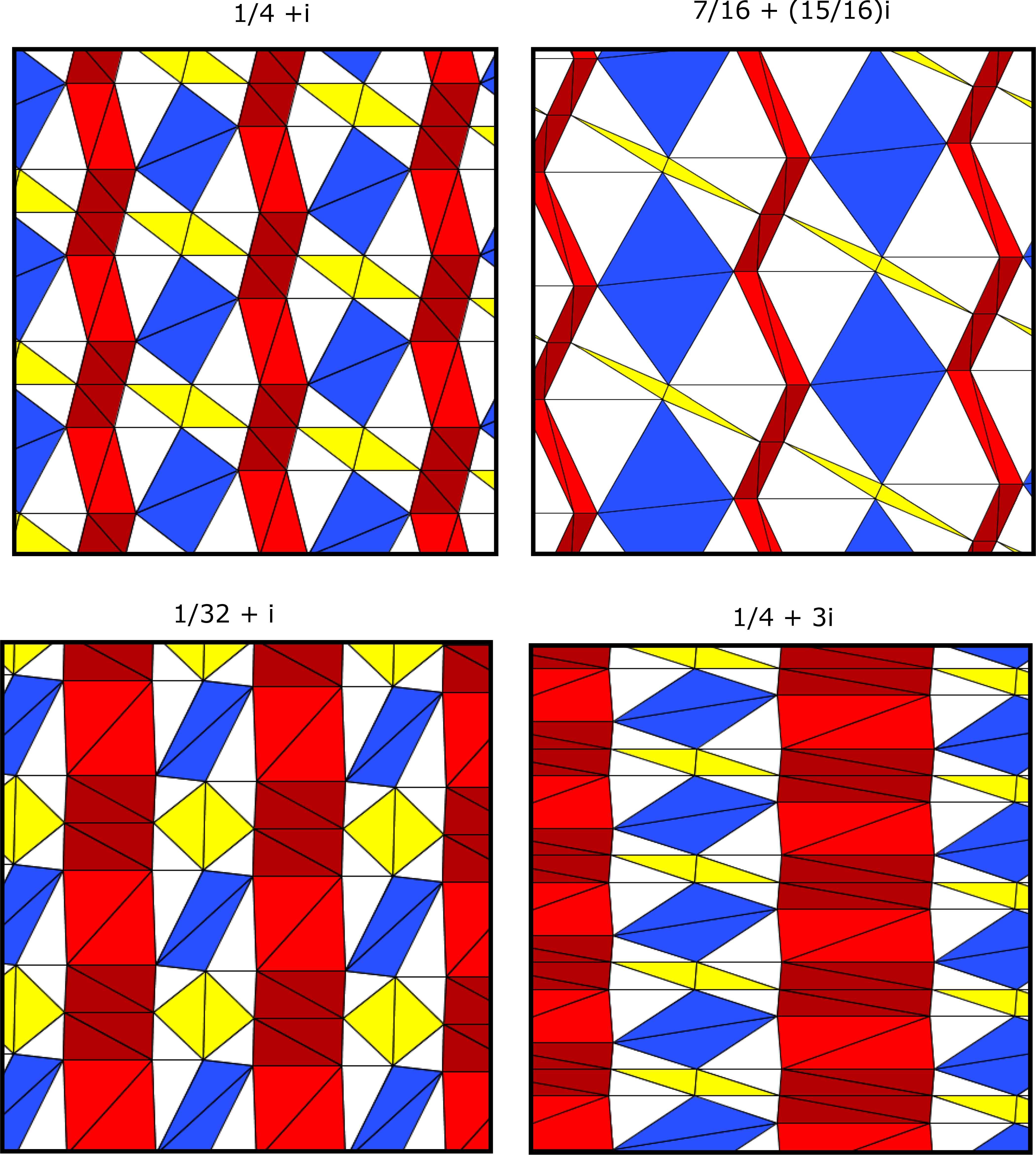}}
\newline
{\bf Figure 2.3:\/}  The tiling at $4$ different parameters
\end{center}

The razor sharp formulas given above
hide a huge amount of experimental
trial and error.  Here we discuss how
we arrived at these equations.
Equation \ref{pup} gives a
truncation of the coordinates
of the embedded pup tent
from [{\bf S\/}].  (An initial version of
[{\bf S\/}] had a slightly different example.)
Only $z_0,z_1,z_2$ are truncated. 
As it is, the object with the exact coordinates
shown in Equation \ref{pup} is
an embedded polyhedral torus
which is within $10^{-31}$ of being flat.

\begin{equation}
  \label{pup}
\begin{matrix}
+0.64 & -0.20 & 0 \\
-1.09 & +0.38 & z_1 \\
-0.25 & +0.51 & z_2 \\
+0.78 & +0.62 & z_3 \\
-0.78 & -0.62 & z_3 \\
+0.25 & -0.51 & z_2 \\
+1.09 & -0.38 & z_1 \\
-0.64 & +0.20 & 0 
\end{matrix}
\hskip 20 pt
\begin{matrix}
\begin{aligned}
z_1 &= 0.0206\ 6632\ 6669\ 8443\ 6159\ 8992\ 3371\ 8861 \\
z_2 &= 0.0048\ 5312\ 7706\ 5192\ 8720\ 4090\ 7479\ 6169 \\
z_3 &= 0.0082\ 2752\ 1455\ 6137\ 1645\ 5791\ 2547\ 8661 
\end{aligned}
\end{matrix}
\end{equation}

The example in Equation \ref{pup}
is only barely embedded.
Some pairs of adjacent triangles in Equation \ref{pup}
are almost completely folded over, giving a very sharp
ridge.  We took such triangle pairs and folded them
completely.  Some other pairs of triangles are very
gently folded.  This happens, for instance, for the ``inner''
$4$ blue triangles in Figure 2.1.  We turned the union of
these $4$ triangles into a paralleogram.

Having degenerated this way, we found the
$2$-parameter family above.   It took us a lot
more trial and error to
find the parametrization above,
given directly in terms of the modular parameter.
We could see from calculations that the vertices of
the intrinsic triangulation ought to have rational
coordinates when suitably translated.
We found these intrinsic coordinates in
Mathematica using a
version of the developing map.

\newpage

  \section{The Ribbon of Highway}

\subsection{Some Positive Functions}

Recall that $\cal I$ is the interior of the bi-cusped fundamental
domain $\cal F$.
Many polynomial expressions in our proofs below
involve the factors
\begin{equation}
  \gamma_0=1-2x, \hskip 12 pt
  \gamma_1=-2x+x^2+y^2, \hskip 12 pt
  \gamma_2=2x-x^2+y^2, \hskip 12 pt
  \gamma_3=2x+x^2+y^2.
\end{equation}
These are all positive on $\cal I$. Compare Equation \ref{DEL}.
Here are two more examples that also appear below.
\begin{equation}
  \gamma_4=2x\gamma_0+(2x+1)(x^2+y^2), \hskip 30 pt
  \gamma_5=2x^{2}\gamma_{0} + x^{2}\gamma_3 +y^{2}\gamma_{3}.
\end{equation}
We will repeatedly use the fact that a positive
polynomial in $x,y, \gamma_0,\gamma_1,\gamma_2,\gamma_3$ is positive
on $\cal I$.

\subsection{Vertical Variations}

\noindent
{\bf Conventions:\/}
Recall that $\rho(u,v,w)=(-u,-v,w)$.
We  work with $\rho$-invariant
$8$ vertex polyhedral
tori, with $\rho$ implementing the vertex permutation
$j \to 7-j$.
Given such an $8$-vertex polyhedral torus $P$,
we let $P_j=(u_j,v_j,w_j)$ denote the $j$th vertex.
Here $j=0,...,7$ as above.  We normalize
so that $w_3=w_4=0$.
\newline

Let
$\theta_j$ denote the cone angle at
vertex $j$, namely the sum of all the
angles, at $P_j$, of the triangles
incident to $P_j$.
By symmetry, $\theta_{7-j}=\theta_j$, and
by the Gauss-Bonnet Theorem,
$\theta_0+\theta_1+\theta_2+\theta_3=8 \pi$.
Thus, the triple $(\theta_0,\theta_1,\theta_2)$
determines all $8$ angles.

We vary $w_0=w_7$ and $w_1=w_6$ and $w_2=w_5$, fixing
all other variables.  Let
\begin{equation}
  \label{ANGLE}
  F(w_0,w_1,w_2)=(\theta_0,\theta_1,\theta_2).
\end{equation}
We let $dF(x,y)$ denote the differential of $F$,
a $3 \times 3$ matrix describing $\partial \theta_i/\partial w_j$,
evaluated at
the points corresponding to the
golden pup tent $P(x+i y)$. 

\begin{lemma}
  $dF(x,y)$ is smooth and invertible everywhere in $\cal I$.
\end{lemma}

\startproof
We compute $dF$ explicitly in
Mathematica and we see that it is a symmetric
matrix whose entries involve rational
functions in $x$ and $y$ and $\sqrt x$.
The denominators factor into powers of $\gamma_3$ and
$\gamma_4$.
We compute
\begin{equation}
\det(dF)
= -\frac{64\sqrt{2}\,x^{3/2} \gamma_5}
{\gamma_3^{4} \gamma_4}.
\end{equation}
Hence $dF$ is smooth and invertible throughout
$\cal I$.
\endproof

\subsection{Proofs modulo the Good Path Lemma}

We use the conventions discussed in the previous section.
We first discuss the three ingredients in the Good Path Lemma,
then state the result.
\newline
\newline
\noindent
{\bf Near Flatness:\/}
Given a polyhedral torus $P$ with cone angles $\{\theta_j\}$  define
\begin{equation}
  \Theta(P)=\max(|\theta_0-2 \pi|, |\theta_1-2 \pi|, |\theta_2-2 \pi|).
\end{equation}
As discussed in the last section,
$\Theta(P)=0$ if and only if all $8$ angles equal $2 \pi$.
We say that $P$ is $\epsilon$-{\it flat\/} if
$\Theta(P) \leq \epsilon$.
\newline
\newline
{\bf Robust Embeddings:\/}
Given polyhedral tori $P$ and $P'$ as above,
we write $P \sim P'$ if these tori
only differ in the
coordinates $w_0,w_1,w_2$ and $w_0',w_1',w_2'$.
Define
\begin{equation}
  \|P-P'\|=\max(|w_0-w_0'|,|w_1-w_1'|,|w_2-w_2'|).
\end{equation}
We say that $P$ is $\lambda$-{\it robustly embedded\/}
if $P'$ is embedded whenever $P \sim P'$ and
$\|P-P'\| \leq \lambda$.
\newline
\newline
{\bf Special Deformations:\/}
Suppose now that $P(z)$ is a
golden pup tent.  We are going to
construct an algebraic map
$(z,t) \to P(z,t)$ which has the property
that $P(z,0)=P(z)$.  The domain of this map
is ${\cal I\/} \times [0,\infty)$.
Really, we are only interested in the image for
very small positive values of $t$.
This map commutes with $\rho$ and is such that
$P_3(z,t)=P_3(z)$ and
$P_4(z,t)=P_4(z)$ for all $t$.
We call any map like this a
{\it special deformation\/}.

\begin{lemma}[Good Path]
  There exists a special deformation
  $P(z,t)$ and positive functions
  $a_z, b_z$ with the following properties.
  \begin{enumerate}
\item For $t>0$ sufficiently small,
    $P(z,t)$ is $a_zt^3$-flat. 
  \item For $t>0$ sufficiently small,
    $P(z,t)$ is $b_zt^2$-robustly embedded.
  \end{enumerate}
 $a_z$ and $b_z$ can be taken independent
  of $z$ if we restrict to a compact subset   of $\cal I$.
\end{lemma}

We fix $z \in \cal I$ and we study polyhedral tori near
$P(z)=P(z,0)$.
For each choice of coordinates $\{u_j\}$ and $\{v_j\}$ we
get a map $F$ as in Equation \ref{ANGLE}.   Since
the map $F$ varies smoothly and $dF$ is invertible on
the golden valley, we see that
$dF$ is invertible in an open neighborhood around
$P(z,0)$.  By the inverse function theorem, the
map $F$ is a local diffeomorphism for any
choice of $\{u_j\}$ and $\{v_j\}$ sufficiently
near the values which lead to $P(z,0)$.

Given the invertibility of the Jacobian,
we can define a flat
polyhedral torus $P'(z,t)$ as long as
$t>0$ is sufficiently small and
$z$ varies within a compact subset of $\cal I$.
The definition is such that $P(z,t) \sim P'(z,t)$.
The Jacobian depends smoothly on all the coordinates, and
the expansion properties of its inverse
likewise vary smoothly.  Combining this fact
with the Good Path Lemma, we get
$\|P(z,t)-P'(z,t)\|<a'_zt^3$ for some $a'_z>0$
which we can take to be locally constant on
compact subsets of $\cal I$.
Given that $P(z,t)$ is $O(t^2)$-robustly embedded,
we see that $P'(z,t)$ is embedded for $t>0$ sufficiently
small.  In short $P'(z,t)$ is a paper torus
which converges to $P(z)$ as $t \to 0$. 

Let $\{{\cal F\/}_n\}$ be a compact exhaustion of $\cal I$
by (say) geodesic triangles, arranged so that
each one is compactly contained in the interior of
the next one.   By compactness and the analysis above there
is some $\epsilon_n>0$ such that $P'(z,t)$ is a paper torus provided
that $z \in {\cal F\/}_n$ and $t \in (0,\epsilon_n)$.
We can choose a continuous positive function
$\chi: {\cal I\/} \to \R$ such that
the restriction of $\chi$ to ${\cal F\/}_n$
is less than $\epsilon_n$.   The union
\begin{equation}
  {\cal U\/} = \bigcup_{z \in {\cal I\/}}\   \bigcup_{t \in (0,\chi(z)]}
  \{P'(z,t)\}.
\end{equation}
is a path connected subset of paper tori.
Here $\{P'(z,t)\}$ denotes the point in $\cal X$
corresponding to $P(z,t)$.
 When we prove the Good Path Lemma below, we
 check that the signs of all the volumes of all the
tetrahedra made from the
vertices of our tori are the same as they are for
the pup tent in Equation \ref{pup}.  Thus $\cal U$
consists of general position pup tents.  In
particular ${\cal U\/} \subset \cal X$.
\newline
\newline
{\bf The Main Theorems:\/}
Define the continuous map
$f_{\chi}:  {\cal I\/} \to \cal M$ as
\begin{equation}
  f_{\chi}(z)=\Phi(P'(z,\chi(z))).
\end{equation}
Here $\cal M$ is the modular surface.
If we make $\chi$ decay fast enough we
can guarantee that $f_{\chi}({\cal I\/}) \subset \cal I$ and that
 $f_{\chi}$ converges to the identity map
on $\partial \cal F$.  But then we have
$\Phi({\cal U\/})=\cal I$.
This proves the Near Universality Theorem.

Given any good polygon $Q$ there is some
$\zeta \in \partial \cal F$ such that the
convex hull of $P(\zeta)$ equals $Q$ up to similarity.
Let $s \mapsto \gamma_s$ be a smooth path
in $\cal I$ from some point we like to $\zeta$.
Up to similarity, 
the path $s \mapsto f_{\chi}(\gamma_s)$
converges to $P(\zeta)$ in the Hausdorff topology as long
as $\chi$ decays fast enough.  This proves the
Collapsibility Theorem.

The Golden Theorem follows immediately from the
definitions and from the proofs we have just given:
Any member of $\cal G$ is the endpoint of a path in
$\cal U$.  Hence $\partial \cal X$ contains a copy of
$\cal G$.  By construction $\Phi({\cal G\/})=\cal M$.

\newpage

\subsection{Proof of the Good Path Lemma}
\label{nearflat}

{\bf Explicit Formulas:\/}
Let $z=x+iy \in \cal I$.
Motivated by experimentation, we
guess the form of our special deformation:
\begin{equation}
\begin{aligned}
u_0(t) &= u_0 + t , & v_0(t) &= v_0 + m t , & w_0(t) &= w_0 + a_0 t^2 \\
u_1(t) &= u_1 + X_1 t^2 , & v_1(t) &= v_1 + X_2 t^2 , & w_1(t) &= w_1 + a_1 t^2 \\
u_2(t) &= u_2 , & v_2(t) &= v_2 + X_1 t^2 , & w_2(t) &= w_2 + a_2 t^2 
\end{aligned}
\end{equation}
We solve $\theta'_j(0)=\theta''_j(0)=0$ for $j=0,1,2$.
We get (and then use)

{\fontsize{10pt}{10pt}\selectfont
\begin{equation}
  m   = \frac{-2 x y}{\gamma_2}, \hskip 30 pt
a_j=
\frac{\alpha_{j} + \alpha_{j1} X_1 + \alpha_{j2} X_2}
{4 \sqrt{2x} \gamma_0 \gamma_1 \gamma_2^2 \gamma_5}.
   \end{equation}
}
{\fontsize{9pt}{9pt}\selectfont
\begin{align*}
\alpha_{0} ={}&
8\,x\,y\,
\bigl(-4x^{2} + 9x^{3} - 7x^{4} - 3x y^{2}  - y^{4}\bigr)\,\gamma_{1},
                \hskip 17 pt
\alpha_{01} =
-4\,y\,\gamma_{0}\,\gamma_{1}\,\gamma_{2}^{2}\,\gamma_{3}^{2}, \hskip
                 17 pt
\alpha_{02} =
2\,(x - 2y^{2})\,\gamma_{0}\,\gamma_{2}^{2}\,\gamma_{3}^{2}, \\[4ex]
\alpha_{1} ={}&
8\,x\,y\,(x^{2} + y^{2})(2x - 3x^{2} + y^{2})\,\gamma_{1}, \hskip 17
                pt
\alpha_{11} ={}
-4\,y\,
\bigl(x^{4} + 6x^{2} + 4x y^{2} + 2x^{2} y^{2} + y^{4}\bigr)\,
\gamma_{0}\,\gamma_{1}\,\gamma_{2}^{2}, \\[1.3ex]
\alpha_{12} ={}&
2\,
\bigl(
2x^{7} \!-\! 9x^{6} \!+\! 12x^{5} \!-\! 4x^{4}
\!+\! 6x^{5} y^{2} \!-\! 11x^{4} y^{2} \!-\! 12x^{3} y^{2} \!-\! 12x^{2} y^{2}
\!+\! 6x^{3} y^{4} \!-\! 3x^{2} y^{4} \!-\! 8x y^{4}
\!+\! 2x y^{6} \!-\! y^{6}
\bigr)\,
\gamma_{0}\,\gamma_{2}^{2}. \\[4ex] 
\alpha_{2} ={}&
4\,x\,y\,
\bigl(-4x^{2} + 6x^{3} - 5x^{4} - 2x y^{2} - 6x^{2} y^{2} - y^{4}\bigr)\,\gamma_{1}, \\[1.3ex]
\alpha_{21} ={}&
4\,x\,
\bigl(2x^{2} - 2x^{3} + x^{4} - 2x y - x^{2} y
      + 2x y^{2} + 2x^{2} y^{2} - y^{3} + y^{4}\bigr)\,
\gamma_{0}\,\gamma_{1}\,\gamma_{2}^{2}, \\[1.3ex]
  \alpha_{22} ={}&
                   \bigl(2x^{3} - x^{4} - 6x y^{2} - 2x^{2} y^{2} - y^{4}\bigr)\,
\gamma_{0}\,\gamma_{2}^{2}\,\gamma_{3}.
\end{align*}
}
Thanks to the the vanishing of the cone angle derivatives, and compactness,
any smooth choice of $X_1,X_2$ leads to a deformation which
satisfies the  flatness condition in the Good Path Lemma.
To get the robust embedding, and motivated by the geometry
we discuss at the end of the proof, we define

{\fontsize{10pt}{10pt}\selectfont
$$
X_1=
\frac{
  -4xy
\left(
\vcenter{
  \hbox{
  $\begin{aligned}[t]
   &40 x^{5} - 60 x^{6} + 30 x^{7} + 25 x^{8} - 15 x^{9}
   - 24 x^{3} y^{2} + 24 x^{4} y^{2} + 50 x^{5} y^{2} + 54 x^{6} y^{2} - 48 x^{7} y^{2} \\
   &+ 20 x^{2} y^{4} + 42 x^{3} y^{4} + 36 x^{4} y^{4} - 54 x^{5} y^{4}
   + 22 x y^{6} + 10 x^{2} y^{6} - 24 x^{3} y^{6} + 3 y^{8} - 3 x y^{8}
  \end{aligned}$}
}\vphantom{\big|}
\right)
}{
  3\gamma_0\gamma_2^2 \gamma_3\Gamma}
  $$
\vspace{-2 pt}
$$
X_2=
\frac{
-4xy \gamma_1\left(
\vcenter{
  \hbox{
  $\begin{aligned}[t]
   &-48 x^{4} + 72 x^{5} - 48 x^{6} - 18 x^{7}
   - 20 x^{5} y - 15 x^{6} y
   - 24 x^{3} y^{2} - 48 x^{4} y^{2} - 30 x^{5} y^{2} \\
   &- 32 x^{2} y^{3} - 40 x^{3} y^{3} - 33 x^{4} y^{3}
   - 6 x^{3} y^{4} - 20 x y^{5} - 21 x^{2} y^{5}
   + 6 x y^{6} - 3 y^{7}
  \end{aligned}$}
}\vphantom{\big|}
\right)
}{
  3\gamma_0\gamma_2^2 \gamma_3\Gamma}
$$
}
\vspace{6 pt}
{\fontsize{9pt}{9pt}\selectfont
$$
\begin{aligned}
\Gamma \;=\;&
  y^{7}
+ 2\gamma_0xy^{6}
+ 8xy^{5} + 7x^{2}y^{5}
+ 16x^{2}y^{4} + 6\gamma_0 x^{3}y^{4} 
+ 11x^{2}y^{3} + \gamma_0 x^{2}y^{3} + 6x^{4}y^{3}
+ 24x^{3}y^{2} \\[2pt]
&+ 24x^{4}y^{2} + 16x^{5}y^{2} + 6\gamma_0 x^{5}y^{2}   + 10\gamma_1 x^{3}y + 5\gamma_1 x^{4}y
+ \tfrac{3}{2}\gamma_0 x^{5} + 6\gamma_0^{2}x^{5} + \tfrac{1}{2}\gamma_0^{3}x^{5} + 6\gamma_0^{2}x^{6} + 12\gamma_0 x^{7}.
\end{aligned}
$$
}

\noindent
{\bf Determinant Calculation:\/}
We let $P=P(x+iy,t)$ be the torus we get with the above choices.
Given  $a<b<c<d\subset \{0,...,7\}$ define
  \begin{equation}
    [abcd] = \det (P_b-P_a,P_c-P_a,P_d-P_a).
  \end{equation}
  This is the $6$ times the signed volume of the tetrahedron defined by these vertices.
  
\begin{lemma}
  \label{ORDER2}
  Let $\F$ be the field of rational functions in $x,y$.
  In all $70$ cases we have
  \begin{equation}\label{poly}
    [abcd] = \sum_{k=0}^6 [abcd]_k\, t^k
        \in \F[t].
  \end{equation}
  All coefficients are finite in $\cal I$,
  the leading coefficient has
  order at most $2$, and the leading coefficient
  does not vanish in $\cal I$.
\end{lemma}

\startproof
In each case we compute that
the expression has the form in Equation \ref{poly}.
Since $[abcd]$ finite on 
${\cal I\/} \times [0,\infty)$, the same goes for
the coefficients. Explicitly, the common denominator of
all the coefficients in all cases is
$\sqrt x \gamma_0 \gamma_1 \gamma_2^3 \gamma_5$.

We compute $45$ and $6$ and $19$ cases respectively where 
{\fontsize{12pt}{12pt}\selectfont
\begin{equation}
  \label{case0}
  [abcd]_0= C \sqrt {x} y^2\gamma_0^{\mu} \gamma_1
  \end{equation}
\begin{equation}
  \label{case1}
 [abcd]_0=0, \hskip 20 pt  [abcd]_1=  \frac{ C\sqrt {x} y^2\gamma_3}{\gamma_1}
\end{equation}
\begin{equation}
  \label{case2}
  [abcd]_0=[abcd]_1=0, \hskip 20 pt
[abcd]_2= \frac{Cx^{3/2} y^2
  \gamma_1^{\mu}}{3\gamma_0\gamma_2\gamma_3^{\mu}}.
\end{equation}
}
Here $\mu \in \{0,1\}$ and  $C/\sqrt 2 \in \Z-\{0\}$.
Only Equation \ref{case2} depends on $X_1,X_2$.
\endproof

\noindent
{\bf Embedding Test:\/} Now we explain how we
test for embeddedness.
In light of Lemma \ref{ORDER2},
it suffices to treat the case when all $70$
determinants are nonzero, which puts the vertices
of $P$ in general position.
  To verify that an edge $(P_a,P_b)$ is disjoint
  from a triangle $(P_c,P_d,P_e)$,  we need to check that
  \begin{equation}
    \label{block}
    \epsilon_1 \sigma_{11} \sigma_{12} + \epsilon_2\sigma_{21}\sigma_{22}
    + \epsilon_3  \sigma_{31}\sigma_{32}<3
    \end{equation}
  for a suitable choice of $\epsilon_1,\epsilon_2,\epsilon_3 \in
  \{-1,1\}$
  and a suitable collection $\{\sigma_{ij}\}$ of signs (which
  are all $\pm 1$) of the expressions $[abcd]$, $[abce]$, $[abde]$,
  $[acde]$, $[bcde]$.
 We call a triple sum in Equation
 \ref{block} a  {\it block\/}.

 \newpage

  Let $\Delta_1$ and $\Delta_2$ be two triangles of $P$.  To test that
  $\Delta_1 \cap \Delta_2=\emptyset$ we need the simultaneous truth of
 $6$ blocks, each saying that a certain edge of one
  triangle
  is disjoint from the other triangle.   To test that
  $\Delta_1 \cap \Delta_2=\{v\}$, a common vertex, we need the
  simultaneous truth of $2$ blocks.
    The general position assumption takes care of pairs of
  triangles that share an edge.
  So, checking that $P$ is embedded
  boils down to checking a clause, the
  {\it embedding clause\/}, made from
$24 \times 6 + 72 \times 2 = 288$ blocks.
The truth of the embedding clause only depends on the
signs of the $70$ determinants.
We call a sign list that satisfies the embedding clause
{\it winning\/}.
\newline
\newline
{\bf Robust Embedding:\/}
Let $\Lambda(z)$ be the list of $70$ signs associated to the
leading coefficients in Lemma \ref{ORDER2}.
We check that $\Lambda(1/4+i)$ is winning, and indeed matches the
sign list for the
embedded pup tent in Equation \ref{pup}.
By Lemma \ref{ORDER2},  the list $\Lambda(z)$ is independent of $z$
and hence winning for all $z \in \cal I$.
Suppose we have  another torus $P^*(z,t)$ with
$\|P^*(z,t)-P(z,t)\|<b_zt^2$.  Given the fact that all
our determinants are smooth functions, Lemma \ref{ORDER2}
tells us that if $b_z$ and $t$ are sufficiently small,
all determinants associated to $P^*(z,t)$ are nonzero, and
the sign list matches $\Lambda(z)$
and hence is winning. So, $P^*(z,t)$ is embedded.  Hence
$P(z,t)$ is $b_z t^2$-robustly
embedded. Given the smooth dependence of everything in sight on $z$,
we can take $b_z$ to be constant on any compact subset
of $\cal I$.
The proof of the Good Path Lemma is done.

\subsection{Discussion and Pictures}

First we discuss the geometric meaning of our
choices for $X_1$ and $X_2$ above.
The $19$ determinants not
independent of $X_1,X_2$  define a union of $7$ lines
on which the order $2$ terms vanish.
One of the triangles cut out by these lines
corresponds to the leading sign list of the pup tent
from Equation \ref{pup}.  We chose
$(X_1,X_2)$ to be its barycenter.  Figure 3.1
shows the lines and the triangle for
parameters $1/4+i$ and $1/8 + (11/8)i$ and $3/8 + (11/8)i$.
In the adjacent triangle, which also would have worked for us.
The sign list is a bit different, but still corresponds to
a pup tent.

\begin{center}
  \resizebox{!}{1.6in}{\includegraphics{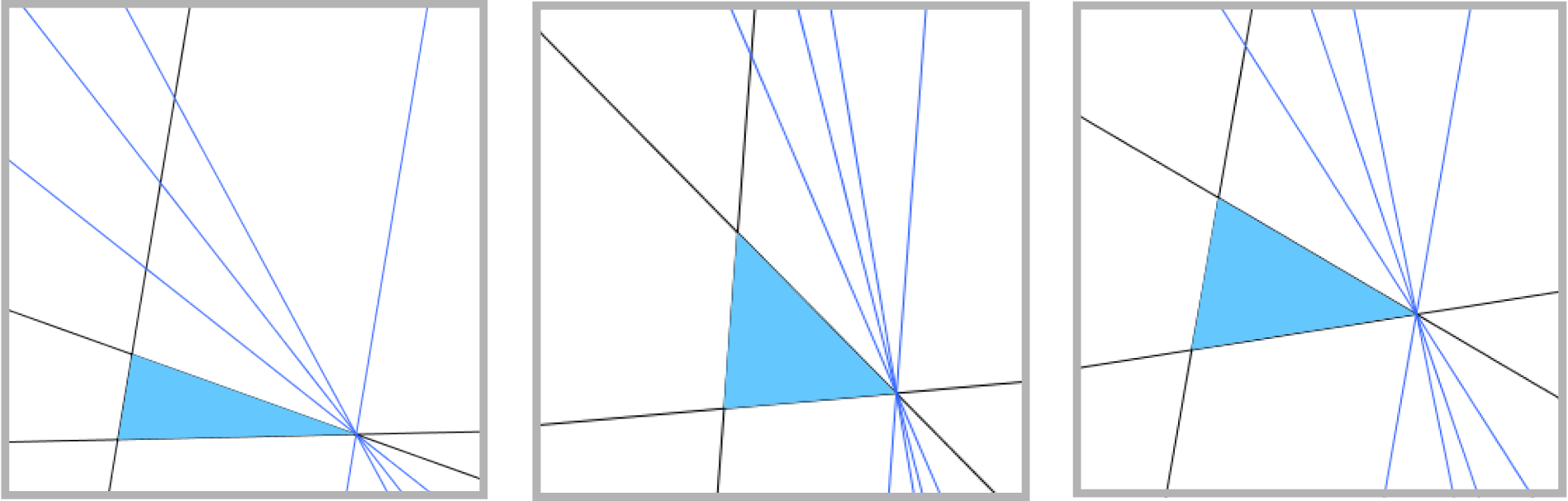}} \newline {\bf
    Figure 3.1:\/} The $7$ lines in the $(X_1,X_2)$-plane and a
  winning blue triangle.
\end{center}

Now we show some pictures of
our constructions in action. In each of the
pictures, we show the
golden pup tent $P(1/4+i,0)$ side by side with the
embedded pup tent $P(1/4+i,1/8)$.
In each picture we are taking  planar slice, and
we are showing the handlebody bounded by the
relevant torus.  The handlebodies respectively
are colored yellow, orange, and red in the
three pictures.

Figure 3.2 shows a slice
by a plane parallel to the $XY$-plane.
There is not really a canonical choice of plane here,
so we just pick one that looks good.
The slice of $P(z,0)$ is on the left and
the slice of $P(z,1/8)$ is on the right.
On the right we have an embedded polygonal disk,
but the disk is just barely embedded.  To make
the point more clear, we show a
close-up of part of the slice.
Our Java program has a zoom feature which
lets you zoom into pictures like this and
see more easily that they are embedded.

\begin{center}
  \resizebox{!}{3.1in}{\includegraphics{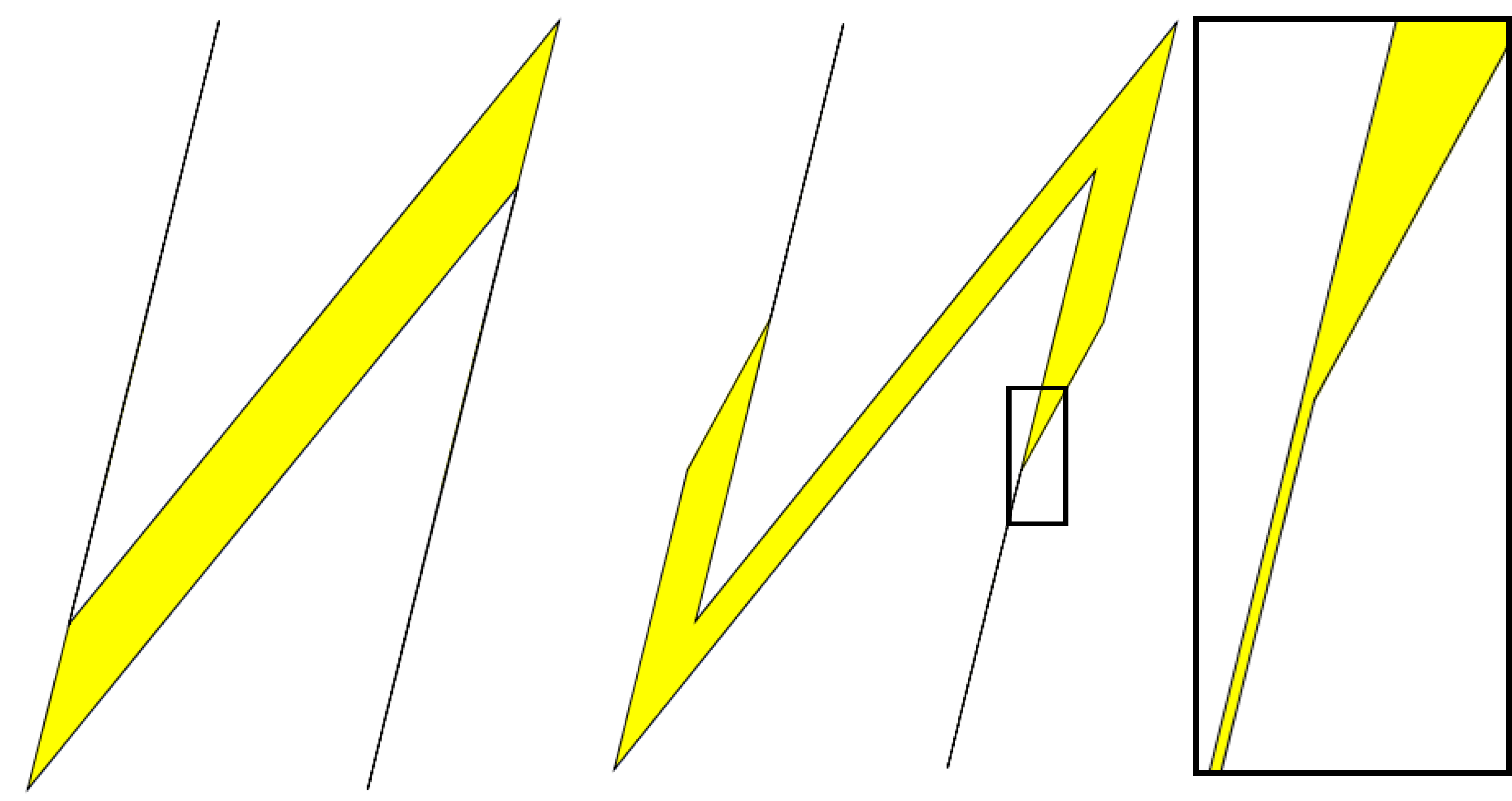}} \newline {\bf
    Figure 3.2:\/} An $XY$-slice for two parameters.
\end{center}

Figure 3.3 has the same format as Figure 3.2.
This time we are slicing by the $XZ$ plane.
This is a canonical choice because the
$XZ$ plane contains the axis of symmetry
of both tori. The figure
in the middle of Figure 3.3
is an embedded annulus.  The two figures on the
right are closeups of the figure in the middle.

\begin{center}
  \resizebox{!}{2.8in}{\includegraphics{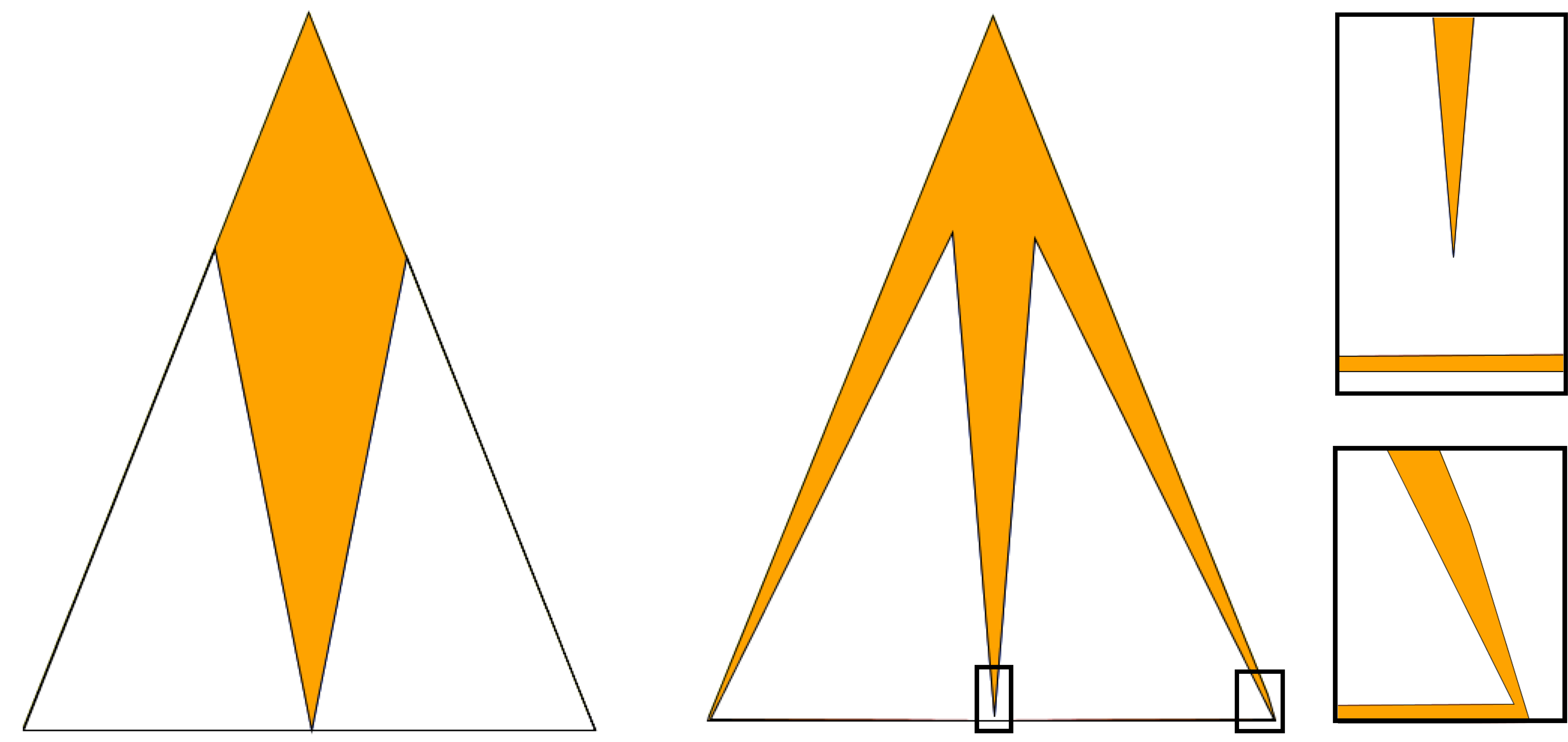}} \newline {\bf
    Figure 3.3:\/} The $XZ$-slice for $P(z,0)$ and $P(z,1/8)$
\end{center}

Figure 3.4 shows the slice by the $YZ$ plane.
Again, this is a canonical choice  because the
$YZ$ plane contains the axis of symmetry.
The figure on the
right is an embedded disk.  This time we are not
showing any closeups. The picture would be similar to
what we see in Figures 3.2 and 3.3.

\begin{center}
  \resizebox{!}{2.2in}{\includegraphics{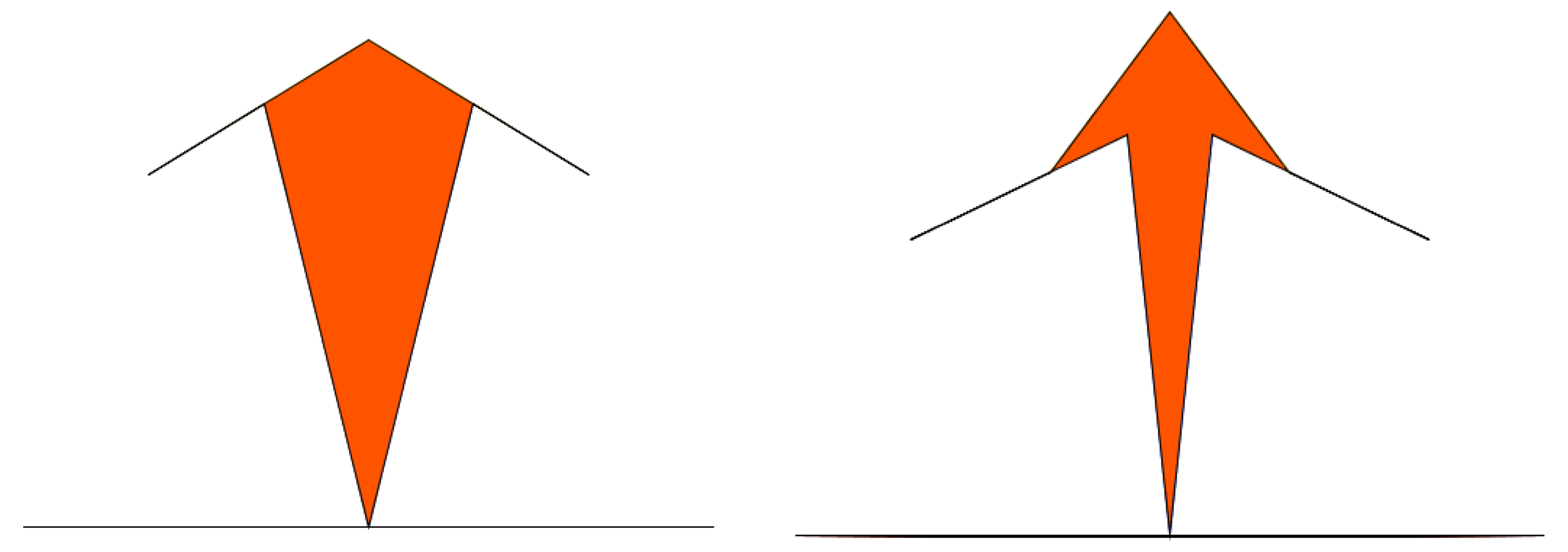}} \newline {\bf
    Figure 3.4:\/} The $YZ$-slice for $P(z,0)$ and $P(z,1/8)$.
\end{center}

Again, we mention that the interested reader can see
and interact with
many more pictures like this using our Java program.

\newpage

  \section{References}

\setlength{\parskip}{0.8\baselineskip}

\noindent
[{\bf ALM\/}] P. Arnoux, S. Leli\`evre, A. M\'alaga, {\it Diplotori:
  a family of polyhedral flat tori.\/} in preparation

\noindent
[{\bf Br\/}]  U. Brehm, Oberwolfach report (1978)

\noindent
[{\bf BZ1\/}] Y. D. Burago and V. A. Zalgaller,
{\it Polyhedral realizations of developments\/} (in Russian)
Vestnik Leningrad Univ. 15, pp 66--80 (1960)

\noindent
[{\bf BZ2\/}] Y. D. Burago and V. A. Zalgaller,
{\it Isometric Embeddings of Two Dimensional Manifolds with a
  polyhedral metric into $\R^3$\/}, Algebra i analiz 7(3) pp 76-95
(1995)  Translation in St. Petersburg Math Journal (3)3, pp 369--385

\noindent
[{\bf LT\/}] F. Lazarus, F. Tallerie, {\it A Universal Triangulation
  for Flat Tori\/}, CS arXiv 2203.05496 (2024)

\noindent
[{\bf Q\/}] T. Quintanar, {\it An explicit embedding of the flat
  square torus in $\E^3$\/}, Journal of Computational Geometry, 11(1):
pp 615--628 (2020)

\noindent
[{\bf S\/}] R. E. Schwartz,  {\it Vertex Minimal Paper Tori\/},
arXiv 2507.14998  (2025)

\noindent
[{\bf S2\/}] R. E. Schwartz,  {\it The Optimal Paper Moebius Band\/},
Annals of Mathematics (2025)

\noindent
[{\bf Se\/}] H. Segerman, {\it Visualizing Mathematics with 3D Printing\/},
Johns Hopkins U. Press  (2016)

\noindent
[{\bf T\/}] T. Tsuboi, {\it On Origami embeddings of flat tori\/},
arXiv 2007.03434 (2020)

\end{document}